\documentclass[11pt]{article}
\usepackage{amsthm, amssymb, amsmath, geometry, verbatim}
\usepackage{authblk}

\numberwithin{equation}{section}
\usepackage{color}
\usepackage{amsfonts,latexsym,amsmath,amscd,geometry,amssymb}
\newcommand \nc{\newcommand}
\newtheorem{theorem}{Theorem}[section]
\newtheorem{lemma}[theorem]{Lemma}

\newtheorem{proposition}[theorem]{Proposition}

\newtheorem{remark}[theorem]{Remark}
\nc{\be}{\begin{eqnarray}}\nc{\ee}{\end{eqnarray}}
\nc{\beq}{\begin{equation}}\nc{\eeq}{\end{equation}}
\nc{\bex}{\begin{eqnarray*}}\nc{\eex}{\end{eqnarray*}}
\nc{\beqno}{\begin{eqnarray*}}\nc{\eeqno}{\end{eqnarray*}}

\newcommand{\R}{\mathbb R}
\newcommand{\dd}{\,\mathrm dx}
\newcommand{\dt}{\,\mathrm dt}
\newcommand{\Div}{\operatorname{div}}
\newcommand{\Curl}{\operatorname{curl}}

\newcommand{\X}{\mathcal X}

\newcommand{\eps}{\varepsilon}
\allowdisplaybreaks

\usepackage[colorlinks,linkcolor=blue]{hyperref}

\begin{document}
\title{Global strong solutions to isentropic
	compressible Navier--Stokes equations
	under a scaling-invariant smallness condition in $\mathbb{R}^2$}
\author{Huanyao Wen\thanks{Email: mahywen@scut.edu.cn}\, and\, Jiayan Wu\thanks{Email: wujiayan@scut.edu.cn}}

\affil{School of Mathematics, South China University of Technology, Guangzhou 510641, China.}

\date{}
\maketitle
\begin{abstract}
In this paper, we consider the Cauchy problem for the isentropic
compressible Navier--Stokes equations with vacuum at far field in two dimensions.
Under an explicit scaling-invariant smallness condition on the initial
data, we establish the global existence and uniqueness of strong solutions. The scaling-invariant smallness condition implies some classes of large data. The main difficulty is to obtain a uniform upper bound for the density under such the condition. Some new estimates of
the dissipation and time-weighted energy estimates with explicit dependence on the initial data are derived, which control the effective viscous flux and close the a priori upper bound for density. In addition, we also obtain algebraic decay of the pressure and the velocity gradient.
	\end{abstract}

 \bigbreak \noindent {\textbf{AMS Subject Classification 2020:}  76N05, 76N10, 35K65.


 \section{Introduction}\label{sec:introduction}
\setcounter{equation}{0}\setcounter{theorem}{0}
\renewcommand{\theequation}{\thesection.\arabic{equation}}
\renewcommand{\thetheorem}{\thesection.\arabic{theorem}}
Compressible Navier--Stokes system is one of the fundamental models arising in fluid dynamics. It describes the motion of viscous, compressible fluids, which attracts a lot of researchers' attention. Much important progress has been made during the past decade. However, there are still some fundamental problems unsolved, such as global existence of smooth solution with large data. This work aims to investigate a related problem for a class of large data, i.e., a scaling-invariant smallness condition of initial quantity. Mathematically, the isentropic compressible Navier--Stokes equations  can be stated as below.
\begin{equation}\label{eq:cns}
	\begin{cases}
		\rho_t+\Div(\rho u)=0,\\
		\rho u_t+\rho u\cdot\nabla u+\nabla P=\mu\Delta u+(\mu+\lambda)\nabla\Div u,\\
	\end{cases}
\end{equation}
for $(x,t)\in\R^d\times[0,\infty)$, where $\rho\ge0$, $u $  denote the density, the
velocity, respectively,  and the pressure is given by  $P(\rho)=A\rho^\gamma$ with $A>0$ and $\gamma>1$.
The viscosity coefficients $\mu$ and $\lambda$  satisfy the following physical conditions:
\begin{equation}\label{eq:parameters}
	\mu>0,\qquad 2\mu+d\lambda\ge0.
\end{equation}	
In the sequel, without loss of generality, we set $A=1$. The system is supplemented with the initial condition
\begin{align}\label{eq:condi-initial}
	(\rho,\rho u)|_{t=0}=(\rho_0,\rho_0u_0),
\end{align}	
and the far-field condition
\begin{align}\label{eq:condi-far}
	(\rho,u)(x,t)\to ( 0, 0) \mbox{ as } |x|\to \infty.
\end{align}

In the absence of vacuum, Serrin~\cite{Serrin1959} and Nash~\cite{Nash1962} established uniqueness result and local wellposedness theory for classical solutions, respectively. Matsumura and Nishida~\cite{MatsumuraNishida1980} obtained global classical solutions for initial data that are small $H^3$ perturbations of a constant equilibrium with positive density. Hoff~\cite{Hoff1995,HoffPolytropic1995} constructed global weak solutions for discontinuous initial data close to a non-vacuum equilibrium. Danchin~\cite{Danchin2000} established global well-posedness for small perturbations in critical Besov spaces. These theories rely on the density being bounded away from zero.  When vacuum is allowed, Lions~\cite{Lions1998} obtained the global finite-energy weak solution for large initial data  provided the exponent $\gamma\geq 9/5$ for three dimensions. See \cite{FeireislEtAl2001, Jiang-Zhang} and \cite{FeireislJDE, B-J2018} for important extensions with respect to the adiabatic index and non-monotone pressure, respectively. The global existence of weak solutions for $\gamma\in(1,\frac{d}{2}]($d=2,3$)$ and the uniqueness of Lions-Feireisl's solution in \cite{Lions1998, FeireislEtAl2001} and Jiang-Zhang's solution in \cite{Jiang-Zhang}. Hu \cite{Hu2,Hu1} investigated the possibilities of the global existence of Lions-Feireisl's solution for $\gamma\in[1,\frac{d}{2}]$. The global existence and uniqueness of the strong solution in three dimensions with small initial data allowing large oscillations was first proved by Huang-Li-Xin \cite{H-L-X:CPAM}. Lei-Xin \cite{Lei-Xin2019} discovered the following scaling invariant property for the first time with respect to the solutions of system (\ref{eq:cns}):
\begin{align}\label{eq:scalling-nature}
	\rho^k(x,t)	=k^{\frac{1}{\gamma}}\rho(k^{\frac{\gamma+1}{2\gamma}}x,kt),\quad u^k(x,t)	=k^{\frac{\gamma-1}{2\gamma}}u(k^{\frac{\gamma+1}{2\gamma}}x,kt).
\end{align}
which means that if $(\rho,u)$ solves the system \eqref{eq:cns}, then $(\rho_k,u_k)$ does, too. The scaling of the equations gives another way to formulate the required smallness. Li~\cite{Li2020Scaling} obtained the global wellposedness of strong solutions under smallness
of a  natural scaling-invariant quantity, with an additional restriction on the
viscosity coefficients $2\mu>\lambda$. Wen~\cite{Wen} removed the additional restriction on the viscosity coefficients in \cite{Li2020Scaling} and for the isentropic system established global strong solutions with far-field vacuum under a small scaling invariant quantity:
$$\bar\rho^3(\|\sqrt{\rho_0}u_0\|_{ 2}^2+\|\rho_0\|_{ \gamma}^\gamma)(\|\nabla u_0\|_{ 2}^2+\|\rho_0^\gamma\|_{ 2}^2)\big[ 1+\bar\rho^{3+\gamma} (\|\sqrt{\rho_0}u_0\|_{ 2}^2+\|\rho_0\|_{ \gamma}^\gamma)\big].$$
See \cite{XuZhong2026MHD, WenXie2026} for further extensions to MHD system and two-fluid model, respectively.

For 2D case, things become quite different, since some useful inequalities are critical in two dimensions. The methods for studying three-dimensional case can not directly be applied to the 2D case. In this context, Luo~\cite{Luo2012} proved local existence of classical solutions to the Cauchy problem of the 2D barotropic compressible Navier--Stokes equations allowing initial vacuum and even compactly supported initial density. Li and Liang~\cite{LiLiang} established a local strong-solution theory based on weighted density norms, finite kinetic energy, and the $L^2$ norm of the velocity gradient. Their main result does not require the usual momentum compatibility condition. 
Li and Xin~\cite{LiXinPre} firstly established global strong solutions and large-time decay under sufficiently small initial total energy. Their argument combines decay of the pressure and velocity gradient with control of the spatial distribution of the mass and the smallness coefficient depends on the initial data.  
Motivated by Lei-Xin's discovery for the scaling-invariant property (\ref{eq:scalling-nature}) for isentropic compressible Navier-Stokes equations, we intend to establish the global well-posedness of the two-dimensional isentropic compressible Navier--Stokes equations with a small initial quantity satisfying (\ref{eq:scalling-nature}).

\bigskip

\subsection*{Notation}
\begin{itemize}
	\item Denote $\displaystyle\int  f {\rm d}x :=\int_{\mathbb{R}^2} f \,{\rm d}x$.
	\item  For $1\le l\le \infty$, we use the following notations for the standard Lebesgue and Sobolev spaces:
	$$\displaystyle L^l=L^l(\mathbb{R}^2), \ \ \ D^{k,l}=\left\{ u\in L^1_{\rm{loc}}(\mathbb{R}^2): \|\nabla^k u \|_{l}<\infty\right\},$$
	$$W^{k,l}=L^l\cap D^{k,l}, \ \ \ H^k=W^{k,2}, \ \ \ D^k=D^{k,2},$$
	$$D_0^1=\Big\{u\in L^6:  \, \|\nabla u\|_{2}<\infty \},$$
	$$\|u\|_{D^{k,l}}=\|\nabla^k u\|_{l}, \mbox{ where } \|\nabla^k u\|_{l}=\|\nabla^k u\|_{L^l}. $$

\item $\dot u=u_t+u\cdot\nabla u.$

  \item $m_0:=\int_{\R^2}\rho_0\dd>0$,\, $P_0:=\rho_0^\gamma$,\, $\bar\rho:=\|\rho_0\|_{\infty}.$

\item $\displaystyle E_0:=\int_{\R^2}\left(\frac12\rho_0|u_0|^2+
		\frac{\rho_0^\gamma}{\gamma-1}\right)\dd.$

\item $\displaystyle H_0:=\|\nabla u_0\|_2^2+\|\rho_0^\gamma\|_2^2.$

\item $\displaystyle R_*:=
		 1+\frac{ \bar\rho^{\frac{\gamma+1}{2}}}{  m_0}\int_{\R^2}\rho_0
		|x| \dd
		+ \sqrt{\frac{2\bar\rho^{1-\gamma}E_0}{m_0}}.$

 \item   Let C denote a generic positive constant depending only on $\gamma,\mu,\lambda$. but independent of $T$ and initial data.
\end{itemize}


We now state the main result.
\begin{theorem}\label{thm:main}
	Let $\eta_0$ be a positive constant and   set
	\begin{equation}
		\bar x=(\mathrm e+|x|^2)^{1/2}\log^{1+\eta_0}(\mathrm e+|x|^2).\label{eq:bar-x}
	\end{equation}
	Assume $\bar x^a\rho_0\in L^1\cap H^1\cap W^{1,q}$ for some $a>1$ and $q>2$, $\sqrt{\rho_0}u_0\in L^2$, and $\nabla u_0\in L^2$. There exists a positive constant
	$\eps_*$ determined by (\ref{mainassum1}) and (\ref{mainassum2}), such that if
	\begin{equation}\label{eq:main-small}
		\begin{aligned}
\bar\rho E_0+\bar\rho E_0\left(1+(\bar\rho E_0)^3\right)(\bar\rho^\gamma m_0)^{-\frac{6}{5}}  R_*^{\frac{4}{5}}\left[(\bar\rho E_0)^2+  [R_*^4+(\bar\rho^\gamma m_0)^2]  (\bar\rho^{1-\gamma}H_0)^2\right]
			\le\eps_*,
		\end{aligned}
	\end{equation}
where $\bar\rho$, $E_0$, $m_0$, $R_*$ and $H_0$ are introduced in Notation, then the
	initial-value problem \eqref{eq:cns}-\eqref{eq:condi-far} admits a unique global strong solution $(\rho, u)$ satisfying for any $0<T<\infty$,
	\begin{equation}\label{eq:main-density}
		0\le\rho(x,t)\le2\bar\rho,
		\qquad x\in\R^2,\quad t\in [0,T],
	\end{equation}
	and
	\begin{equation}
		\left\{
		\begin{aligned}\label{eq:main-reg}
			&\rho\in C([0,T];L^1\cap H^1\cap W^{1,q}),
			\qquad \bar x^a\rho\in L^\infty(0,T;L^1\cap H^1\cap W^{1,q}),
			\\
			&  \nabla u, \sqrt\rho u, \bar x^{-1} u, \sqrt t\sqrt\rho\dot u\in L^\infty(0,T;L^2),\\
			&\nabla u\in  L^2(0,T;H^{1})\cap L^{\frac{q+1}{q}}(0,T;W^{1,q}), 	\qquad \sqrt t \nabla u\in  L^2(0,T;W^{1,q}),\\
			& \sqrt\rho \dot u, \sqrt t  \bar x^{-1}  \dot u, \sqrt t\nabla\dot u \in L^2(0,T;L^2).
		\end{aligned}
		\right.
	\end{equation}
	Moreover, for $t\ge1$, the solution satisfies
	\begin{equation}
		\|\nabla u\|_2+\|P\|_2\le\frac{C\sqrt{E_0}}{t^{1/2}},\qquad
		\|\nabla F\|_2+\|\nabla\omega\|_2\le\frac{C\sqrt{\bar\rho E_0}}{t},
	\end{equation}
	where
	$F=(2\mu+\lambda)\Div u-P$ and
	$\omega=\partial_1u^2-\partial_2u^1$.
\end{theorem}

\begin{remark}
	It is easy to verify that   $\bar\rho E_0$, $\bar\rho^{1-\gamma}H_0$,
	$\bar\rho^\gamma m_0$, and $R_*$  are invariant under the scaling \eqref{eq:scalling-nature}. Therefore the small quantity  in (\ref{eq:main-small}) is scaling invariant. When this work was finished, the authors had a friendly private discussion with Professor Boqiang L\"u and found that, under another interesting scaling-invariant smallness condition regarding the following quantity:
 \begin{equation*}
\left[\bar\rho^\gamma m_0+(\bar\rho^\gamma m_0)^{-1}+\bar\rho^{\gamma+\frac{a(\gamma+1)}{2}}\|\rho_0|x|^a\|_{L^1}\right](1+\bar\rho E_0)^\frac{3}{2}\bar\rho^{1-\gamma}\left(\|\nabla u_0\|_{L^2}^2+\|\rho_0^\gamma\|_{L^2}^2\right),
 \end{equation*}
 and using different techniques, they had independently obtained results on the global existence and uniqueness of strong solutions to system (\ref{eq:cns}) (\cite{Lu-Song}). Furthermore, neither of these two conditions implies the other.
\end{remark}

\begin{remark}\label{rem:compatibility}
	If	the initial data are further assumed to satisfy $ \sqrt{\rho_0}u_0\in L^2,
	\nabla^2 u_0\in L^2$ and the following compatibility condition
	\begin{align}
		\mu\Delta u_0+(\mu+\lambda)\nabla\Div u_0-\nabla P_0
		=\sqrt{\rho_0}\,g_0,\qquad
		g_0\in L^2,\label{eq:compat}
	\end{align}
	we can get higher regularity of the solution:
	\begin{equation}
		\begin{aligned}\label{eq:main-reg-strong}
			\nabla^2 u, \sqrt\rho\dot u \in L^\infty(0,T;L^2),
			\nabla\dot u \in L^2(0,T;L^2).
		\end{aligned}
	\end{equation}
\end{remark}

We briefly present our main ideas in the proof of the main result.  Under the assumptions on the initial data stated in Theorem \ref{thm:main}, the local existence and uniqueness of the solution to the Cauchy problem \eqref{eq:cns}-\eqref{eq:condi-far} have been established. To prove the global existence, it suffices to derive some a priori estimates globally in time. The key estimate is the upper bound for density. One of the main challenges compared with the three-dimensional case (\cite{Wen}) is that some inequalities are critical. Therefore some weighted in space and in time estimated are necessary to develop under the scaling-invariant smallness condition (\ref{eq:main-small}), see Lemmas \ref{lem:riccati}, \ref{lem:second-energy} and \ref{eq:increaseR}. With those estimates together with the following crucial estimate (see Lemma \ref{lem:bridge}):
$$
		\begin{aligned}
		\int_{\R^2}\rho|v|^p\dd
		\leq C\bar\rho^{-\frac{\gamma}{5} }
		\left(m_0R\right)^{\frac{4}{5}}\left( m_0^{-1/2}\|\sqrt\rho v\|_2
		+[(\bar\rho^\gamma m_0)^{-1/2}R+1]\|\nabla v\|_2\right)^p.
	\end{aligned}
$$ for any given $v\in D^{1,2}$ and $\sqrt{\rho}v\in L^2$, we apply Zlotnik's inequality to close the a priori bound for density from $2\bar\rho$ to $\frac{7\bar\rho }{4}$, see Proposition \ref{prop:density}.

The rest of the paper is organized as follows.
In Section~\ref{sec:tools}, we present some preliminary results which will be used in derivation of global a priori estimates.
In Section~\ref{sec:proof}, we establish a series of global-in-time energy estimates and in particular obtain the crucial upper bound for density, and then finish the proof of the main result.

}

\section{Preliminaries}\label{sec:tools}
	In this section,   we will recall some known facts and elementary inequalities which will be used frequently later.




\begin{lemma}[\cite{GN-ineq}] \label{GN-inequality}
	Let $1<m<\infty$ and $1\le r<\infty$. For $f\in L^r(\R^2)$ with
	$\nabla f\in L^m(\R^2)$,
	\[
	\|f\|_q\le C\|\nabla f\|_m^\theta\|f\|_r^{1-\theta},
	\qquad
	\theta=\frac{1/r-1/q}{1/r-1/m+1/2},
	\]
	in the following ranges: if $m<2$, assume $r<2m/(2-m)$ and
	$r\le q\le2m/(2-m)$; if $m=2$, assume $r\le q<\infty$;
	if $m>2$, assume $r\le q\le\infty$. The constant $C$ depends on $m,r,q$.
\end{lemma}

Taking divergence and curl in the momentum equation gives
\begin{equation}\label{eq:elliptic}
	\Delta F=\Div(\rho\dot u),\qquad
	\mu\Delta\omega=\Curl(\rho\dot u),
\end{equation}
where
\[
F=\nu\Div u-P,\qquad
\omega=\partial_1u^2-\partial_2u^1.
\]
Then the standard $L^p$-estimate for the elliptic system \eqref{eq:elliptic} yields the following estimate.
\begin{lemma}\emph{\cite[Lemma 2.5]{LiXinPre}}\label{lem:elliptic}
	For $2\leq p<\infty$, there holds that
	\begin{equation}\label{eq:elliptic-Lp}
		\|\nabla F\|_p+\mu\|\nabla\omega\|_p\le C_p\|\rho\dot u\|_p
	\end{equation}
	and that
	\begin{equation}\label{eq:flux-Lp}
		\|F\|_p+\|\omega\|_p
		\le C_p\|\rho\dot u\|_2^{1-2/p}H^{1/p}.
	\end{equation}
	If $0\le\rho\le2\bar\rho$, then
	\begin{equation}\label{eq:velocity-Lp}
		\|\nabla u\|_p
		\le C_p\bar\rho^{1/2-1/p}H^{1/p}
		\|\sqrt\rho\dot u\|_2^{1-2/p}+C_p\|P\|_p,
	\end{equation}
	and, in particular,
	\begin{equation}\label{eq:grad-L3}
		\|\nabla u\|_3^3
		\le C\bar\rho^{1/2}H\|\sqrt\rho\dot u\|_2+C\|P\|_3^3,
	\end{equation}
	\begin{equation}\label{eq:grad-L4}
		\|\nabla u\|_4^4
		\le C\bar\rho H\|\sqrt\rho\dot u\|_2^2+C\|P\|_4^4.
	\end{equation}
		Here $H =\|\nabla u \|_2^2+\|P \|_2^2$ and $C_p$ in \eqref{eq:elliptic-Lp} depends only on $p$, and constants
	in \eqref{eq:flux-Lp}--\eqref{eq:grad-L4} may also depend on $\mu,\lambda$.
\end{lemma}

\begin{lemma}\emph{\cite[Theorem B.1]{Lions-book1}}\label{lem:weighted-sobolev}
	For $m\in [2,\infty)$ and $\theta \in (1 + m/2,\infty)$, there exists a positive constant
	$C_m$ such that for all $w\in D^1(\mathbb R^2)$,
	\begin{equation}\label{eq:weighted-sobolev}
		\left(\int_{\R^2} \frac{|w|^m}{\mathrm e+|x|^{2}}[\log (\mathrm e+|x|^2)]^\theta\dd\right)^{1/m}
		\le C_m\bigl(\|w\|_{L^2(B_1)}+\|\nabla w\|_2\bigr),
		\qquad w\in D^{1,2}(\R^2).
	\end{equation}
\end{lemma}

\medskip

Finally, the following Zlotnik's  inequality  plays an important role to the uniform upper bound of the density.
\begin{lemma}[\cite{Zlotnik}] \label{lem:zlotnik}
	Let the function $y$ satisfy
	\[
	y'(t)=g(y(t))+b'(t)
	\quad \text{on } [0,T],
	\qquad y(0)=y^0,
	\]
	where $g\in C(\mathbb{R})$ and $y,b\in W^{1,1}(0,T)$.
	Suppose that $\lim_{\xi\to+\infty}g(\xi)=-\infty$
	and
	\begin{equation}\label{eq:zlotnik-increment}
		b(t_2)-b(t_1)
		\le N_0+N_1(t_2-t_1)
	\end{equation}
	for all $0\le t_1<t_2\le T$, with some constants
	$N_0\ge0$ and $N_1\ge0$. Then
	\begin{equation}\label{eq:zlotnik-bound}
		y(t)\le \max\{y^0,\bar{\zeta}\}+N_0<\infty
		\quad \text{on } [0,T],
	\end{equation}
	where $\bar{\zeta}$ is a constant such that
	\[
	g(\zeta)\le -N_1
	\qquad \text{for all } \zeta\ge\bar{\zeta}.
	\]
\end{lemma}

	\section{Proof of Theorem~\ref{thm:main}}\label{sec:proof}
	 We first recall a local existence theorem.
	\begin{proposition}\emph{\cite[Theorem 1.1]{LiLiang}}\label{lem:local}
		Assume that the coefficients satisfy \eqref{eq:parameters}.
		Let
		\[
		0<\eta_0\le1,\qquad q>2,\qquad 1<a<2,
		\]
		and define
		\[
		\bar x
		:=
		\bigl(\mathrm e+|x|^2\bigr)^{1/2}
		\log^{1+\eta_0}\bigl(\mathrm e+|x|^2\bigr).
		\]
		Suppose that the initial data satisfy
		\begin{equation}\label{eq:cited-unit-mass}
			\int_{\mathbb{R}^2}\rho_0(x)\,\mathrm dx=1
		\end{equation}
		and
		\begin{equation}\label{eq:cited-local-data}
			\begin{aligned}
				&\rho_0\ge0,\qquad
				\bar x^a\rho_0\in L^1\cap H^1\cap W^{1,q},\\
				&\nabla u_0\in L^2,\qquad
				\sqrt{\rho_0}\,u_0\in L^2.
			\end{aligned}
		\end{equation}
		Then there exists $T_0>0$ such that the Cauchy problem
		for \eqref{eq:cns} admits a unique local strong solution
		$(\rho,u)$ in $\mathbb{R}^2\times(0,T_0]$, satisfying
		\begin{equation}\label{eq:cited-local-regularity}
			\left\{ \begin{aligned}
				&\rho\in C\bigl([0,T_0];L^1\cap H^1\cap W^{1,q}\bigr),\quad
				\bar x^a\rho
				\in L^\infty\bigl(0,T_0;L^1\cap H^1\cap W^{1,q}\bigr),\\
				&\sqrt{\rho}\,u,\quad
				\nabla u,\quad
				\bar x^{-1}u,\quad
				\sqrt{t}\sqrt{\rho}\,u_t
				\in L^\infty(0,T_0;L^2),\\
				&\nabla u
				\in L^2(0,T_0;H^1)
				\cap L^{(q+1)/q}(0,T_0;W^{1,q}),\quad
				\sqrt{t}\nabla u
				\in L^2(0,T_0;W^{1,q}),\\
				&\sqrt{\rho}\,u_t,\quad
				\sqrt{t}\nabla u_t,\quad
				\sqrt{t}\bar x^{-1}u_t
				\in L^2(0,T_0;L^2).
			\end{aligned}
			\right.
		\end{equation}
		Moreover, there exists a constant $N>0$ such that
		\begin{equation}\label{eq:cited-local-mass}
			\inf_{0\le t\le T_0}
			\int_{B_N}\rho(x,t)\,\mathrm dx
			\ge\frac14,
		\end{equation}
		where $B_N=\{x\in\mathbb{R}^2:|x|<N\}$.
	\end{proposition}
	
With the local existence result, the maximal existence time of the strong
solution $T^*$ must be positive. One of the key estimates is to derive the upper bound of density, namely,
\begin{equation}\label{eq:bootstrap}
	0\le\rho(x,t)\le2\bar\rho,
\end{equation} for all $x\in\R^2$ and $t\in[0, T^*)$.

The proof of (\ref{eq:bootstrap}) consists of the following two propositions. The first one gives the energy,
dissipation, and decay bounds under \eqref{eq:bootstrap} on some time interval $[0,T]$; the second
improves the density bound to $7\bar\rho/4$ on the same time interval. Then (\ref{eq:bootstrap}) follows from the standard continuity arguments.

\subsection{Energy and dissipation estimates}\label{sec:energy}


\begin{proposition}\label{prop:energy}
	Letting $(\rho,u)$ be the strong solution on $[0,T]$, there exists a positive constant $\eps_1$ depending only on $\mu,\lambda,\gamma$ and some other known parameters but independent of T and the initial data, such that if
	\[
	0\le\rho(x,t)\le2\bar\rho,\qquad x\in\R^2,\quad 0\le t\le T,
	\]
	then there holds that
\begin{equation}\label{eq:summary-one}
		\sup_{0\le t\le T}H(t)
		+\int_0^T\left(\|\sqrt\rho\,\dot u\|_2^2+\|P\|_3^3\right)\dt
		\le CH_0,
	\end{equation}
	\begin{equation}\label{eq:summary-decay}
		H(t)\le\frac{CE_0}{t+E_0/H_0},\qquad
		\|\sqrt\rho\,\dot u(t)\|_2^2+\|P(t)\|_3^3
		\le\frac{CE_0}{t(t+E_0/H_0)},\qquad 0<t\le T,
	\end{equation}
	\begin{equation}\label{eq:summary-two}
		\sup_{0<t\le T}t\left(\|\sqrt\rho\,\dot u \|_2^2+\|P \|_3^3\right)
		+\int_0^Tt\left(\|\nabla\dot u\|_2^2+\|P\|_4^4\right)\dt
		\le CH_0,
	\end{equation}
	and
	\begin{equation}\label{eq:summary-three}
		\begin{aligned}
			\int_0^T H(t)\dt
			+\int_0^Tt\left(\|\sqrt\rho\,\dot u\|_2^2+\|P\|_3^3\right)\dt
			+\int_0^Tt^2\left(\|\nabla\dot u\|_2^2+\|P\|_4^4\right)\dt
			\le C(E_0+m_0H_0),
		\end{aligned}
	\end{equation}
	where $H=\|\nabla u \|_2^2+\|P \|_2^2$, provided that
\begin{equation}\label{mainassum1}\bar\rho E_0 \le \eps_1,\end{equation} and that $\eps_1$ is determined by (\ref{assum1}).
\end{proposition}

\begin{lemma} \label{lem:basic}
	The solution satisfies
	\begin{equation}\label{eq:rho}
		m:=\int\rho(x,t)\dd=m_0
	\end{equation}
	and
	\begin{equation}\label{eq:basic}
		E(t)+\int_0^t
		\left(\mu\|\nabla u\|_2^2
		+(\mu+\lambda)\|\Div u\|_2^2\right)\,\mathrm ds=E_0,
	\end{equation}
	where
	$\displaystyle E(t)=\int_{\R^2}\left(\frac12\rho|u|^2+\frac{P}{\gamma-1}\right)\dd$.
	
\end{lemma}

The following pressure identity will be used in the energy estimates.
\begin{lemma} \label{lem:pressure}
	For every fixed $p\ge2$, we have
	\begin{equation}\label{eq:pressure-power}
		\frac{\mathrm d}{\mathrm dt}\|P\|_p^p+
		\frac{p\gamma-1}{2\mu+\lambda}\|P\|_{p+1}^{p+1}
		=-\frac{p\gamma-1}{2\mu+\lambda}\int P^pF\dd.
	\end{equation}
\end{lemma}
\begin{proof}
	The pressure satisfies
	\begin{equation}\label{P-iden}(\partial_t+u\cdot\nabla)P=-\gamma P\Div u.\end{equation}  Multiplying \eqref{P-iden} by $pP^{p-1}$, we get
	\begin{equation}\label{P-eq1}
	\partial_tP^p+\Div(uP^p)=-(p\gamma-1)P^p\Div u.
	\end{equation}
	Integrating \eqref{P-eq1} in $x$ and substituting $\Div u=(F+P)/\nu$ proves
	\eqref{eq:pressure-power}.
\end{proof}

\begin{lemma} \label{lem:first-energy}
	Under the assumption of Proposition \ref{prop:energy}, we have
	\begin{equation}\label{eq:first-energy}
		H(\beta)+\int_\alpha^\beta
		\left(\|\sqrt\rho\,\dot u\|_2^2+\|P\|_3^3\right)\dt
		\le CH(\alpha),\qquad 0\le\alpha<\beta\le T,
	\end{equation}
	provided that
	\begin{equation*}\bar\rho E_0 \le \min\left\{1,  \frac{c}{2C}\right\}.\end{equation*}
	In  particular, we have
	\begin{equation}\label{eq:first-integrated}
		\sup_{0\le t\le T}H(t)
		+\int_0^T\left(\|\sqrt\rho\,\dot u\|_2^2+\|P\|_3^3\right)\dt
		\le CH_0.
	\end{equation}
\end{lemma}
\begin{proof}
	By a calculation similar to that in \cite{LiXinPre}, using \eqref{eq:grad-L3} and Young's inequatity, we obtain
	\begin{equation}\label{eq:B-rough}
		\begin{aligned}
			&\frac{\mathrm d}{\mathrm dt}
			\left[\frac\mu2\|\nabla u\|_2^2
			+\frac{\mu+\lambda}{2}\|\Div u\|_2^2
			-\int P\Div u\dd\right]
			+\|\sqrt\rho\,\dot u\|_2^2\\
			&\qquad\le C\|\nabla u\|_{3}^3+C\int P|\nabla u|^2\dd\\
			&\qquad\le \frac12\|\sqrt\rho\,\dot u\|_2^2
			+C\bar\rho H^2+C\|P\|_3^3.
		\end{aligned}
	\end{equation}
	Taking $p=2$ in \eqref{eq:pressure-power} and using
	\eqref{eq:flux-Lp} and $P^2|F|\le\eta P^3+C_\eta|F|^3$, with $\eta>0$ sufficiently small,
	we obtain, for any $\delta>0$,
	\begin{equation}\label{eq:V-energy}
		\begin{aligned}
			\frac{\mathrm d}{\mathrm dt}\|P\|_2^2+\frac{2\gamma-1}{2(2\mu+\lambda)}\|P\|_3^3&\leq  C \|F\|_3^3\\
			&\leq\delta\|\sqrt\rho\,\dot u\|_2^2+C_\delta\bar\rho H^2.
		\end{aligned}
	\end{equation}
	Adding $K_2$ times \eqref{eq:V-energy} to \eqref{eq:B-rough}, we get
	\[
	\frac{\mathrm d}{\mathrm dt}\X
	+\frac12\|\sqrt\rho\,\dot u\|_2^2+K_2\frac{2\gamma-1}{2(2\mu+\lambda)}\|P\|_3^3
	\le K_2\delta\|\sqrt\rho\,\dot u\|_2^2
	+C_{\delta}\bar\rho H^2+C\|P\|_3^3,
	\]
	where
	\begin{equation}\label{eq:X-definition}
		\X(t)=\frac\mu2\|\nabla u\|_2^2
		+\frac{\mu+\lambda}{2}\|\Div u\|_2^2
		-\int P\Div u\dd+K_2\|P\|_2^2.
	\end{equation}
	By Young's inequality,
	\[
	\left|\int P\Div u\dd\right|
	\le\frac\mu4\|\nabla u\|_2^2+C\|P\|_2^2,
	\]
	we first choose $K_2$ sufficiently large such that $cH\le\X\le CH$ and
	$K_2\frac{2\gamma-1}{2\nu}\ge C+1$. We then fix $\delta$ such that
	$K_2\delta\le1/4$.
	It follows that
	\begin{equation}\label{eq:X-before-absorb}
		\frac{\mathrm d}{\mathrm dt}\X
		+c\left(\|\sqrt\rho\,\dot u\|_2^2+\|P\|_3^3\right)
		\le C\bar\rho H^2,
	\end{equation}
	where the constants are independent of the initial data and $T$.
	Next,  by using the interpolation inequality
	 and \eqref{eq:basic} imply
	\begin{equation}\label{eq:pressure-key}
		\begin{aligned}
			H^2
			&\le 2\|\nabla u\|_2^2H+\|P\|_2^4\\
			&\le C\|\nabla u\|_2^2\X+\|P\|_1\|P\|_3^3\\
			&\le C\|\nabla u\|_2^2\X+CE_0\|P\|_3^3.
		\end{aligned}
	\end{equation}
	Substituting \eqref{eq:pressure-key} into
	\eqref{eq:X-before-absorb}, we obtain
	\begin{equation}
		\begin{aligned}
			\frac{\mathrm d}{\mathrm dt}\X
			+c\left(\|\sqrt\rho\,\dot u\|_2^2+\|P\|_3^3\right)
			&\leq C\bar\rho\|\nabla u\|_2^2\X
			+C\bar\rho E_0\|P\|_3^3,
		\end{aligned}
	\end{equation}
    which gives
    \begin{equation}\label{eq:X-coercive}
    	\frac{\mathrm d}{\mathrm dt}\X
    	+c\left(\|\sqrt\rho\,\dot u\|_2^2+\|P\|_3^3\right)
    	\le C\bar\rho\|\nabla u\|_2^2\X,
    \end{equation}
      provided that
      \begin{equation}\label{assum1}\bar\rho E_0 \le \min\left\{1,  \frac{c}{2C}\right\}.\end{equation}
	
	For $0\le\alpha\le t\le T$, define
	\begin{equation}\label{M}
	M_\alpha(t)
	=\exp\left(-C\bar\rho\int_\alpha^t\|\nabla u(s)\|_2^2\,\mathrm ds\right).
	\end{equation}
	By \eqref{eq:basic} and $\bar\rho E_0\le1$,
	$e^{-C/\mu}\le M_\alpha(t)\le1$. Hence
	\begin{equation}\label{M-X}
	\frac{\mathrm d}{\mathrm dt}(M_\alpha\X)
	+cM_\alpha\left(\|\sqrt\rho\,\dot u\|_2^2+\|P\|_3^3\right)\le0.
	\end{equation}
	Integrating \eqref{M-X} over $[\alpha,\beta]$, and using $cH\le\X\le CH$, we get
	\[
	H(\beta)+\int_\alpha^\beta
	\left(\|\sqrt\rho\,\dot u\|_2^2+\|P\|_3^3\right)\dt
	\le CH(\alpha).
	\]
	Then letting $\alpha\to 0^+$,    we obtain \eqref{eq:first-integrated} by taking the supremum
	over $\beta\in[0,T]$.  This completes the proof. \end{proof}

\begin{lemma}\label{lem:riccati}
	Under the assumption of Proposition \ref{prop:energy}, we have
	\begin{equation}\label{eq:H-decay}
		H(t)\le\frac{CH_0}{1+tH_0/E_0},
	\end{equation} provided that
	\begin{equation*}\bar\rho E_0 \le \min\left\{1,  \frac{c}{2C}\right\}.\end{equation*}
\end{lemma}
\begin{proof}
	By \eqref{eq:pressure-key} and $cH\le\X\le CH$, we have
	\begin{equation}\label{P-es}
	\X^2\le C\|\nabla u\|_2^2\X+CE_0\|P\|_3^3.
	\end{equation}
	Using \eqref{eq:X-coercive}
 and  \eqref{P-es}, we obtain
	\begin{equation}\label{differential-inequality}
		\frac{\mathrm d}{\mathrm dt}\X+\frac c{E_0}\X^2
		\le C(\bar\rho+E_0^{-1})\|\nabla u\|_2^2\X.
	\end{equation}
	Set
	\[
	Q(t)=C(\bar\rho+E_0^{-1})
	\int_0^t\|\nabla u(s)\|_2^2\,\mathrm ds.
	\]
	By \eqref{eq:basic} and $\bar\rho E_0\le1$, then
	$0\le Q(t)\le C(\bar\rho E_0+1)/\mu\le C$.
	Multiplying  \eqref{differential-inequality} by $e^{-Q}$ gives
	\[
	\frac{\mathrm d}{\mathrm dt}(e^{-Q}\X)
	+\frac c{E_0}e^Q(e^{-Q}\X)^2\le0.
	\]
	Since $e^Q\ge1$, integrating the above inequality over
	$[\alpha,t]$ for $0<\alpha<t$, yields
	\[
	e^{-Q(t)}\X(t)
	\le
	\frac{e^{-Q(\alpha)}\X(\alpha)}
	{1+c(t-\alpha)e^{-Q(\alpha)}\X(\alpha)/E_0}.
	\]
	By Lemma~\ref{lem:first-energy}, $\X(\alpha)\le CH_0$ uniformly in $\alpha$.
	Since the function $$z\mapsto\frac{ z}{(1+c(t-\alpha)z/E_0)}$$ is increasing, it follows that
	\[
	e^{-Q(t)}\X(t)
	\le\frac{CH_0}{1+ctH_0/E_0}
	\]
	after letting $\alpha\downarrow0$.
	Finally, using the uniform bound for $Q$ and $cH\le\X$, we obtain  \eqref{eq:H-decay}.
\end{proof}

\begin{lemma}\label{lem:second-energy}
	Under the assumption of Proposition \ref{prop:energy}, we have
	\begin{equation}\label{eq:Y-t}
		\sup_{0<t\le T}t[\|\sqrt\rho\,\dot u(t)\|_2^2+ \|P(t)\|_3^3]
		+\int_0^Tt\left(\|\nabla\dot u\|_2^2+\|P\|_4^4\right)\dt
		\le CH_0,
	\end{equation}
	 provided that
	\begin{equation*}\bar\rho E_0 \le \min\left\{1,  \frac{c}{2C}\right\}.\end{equation*}
	Moreover,
	\begin{equation}\label{eq:Y-decay}
		\|\sqrt\rho\,\dot u(t)\|_2^2+ \|P(t)\|_3^3\le\frac{CE_0}{t(t+H_0/E_0)},\qquad
		\sup_{0<t\le T}t^2[\|\sqrt\rho\,\dot u(t)\|_2^2+ \|P(t)\|_3^3]\le CE_0,
	\end{equation}
\end{lemma}
\begin{proof}
	By a calculation similar to that in \cite{LiXinPre},  using \eqref{eq:grad-L4} and Young's inequatity,  we obtain
	\begin{equation}\label{eq:accestimate}
		\begin{aligned}
			\frac{\mathrm d}{\mathrm dt}\|\sqrt\rho\,\dot u\|_2^2
			+\mu\|\nabla\dot u\|_2^2&\leq C\|\nabla u\|_4^4+C\|P\|_4^4\\
			&\leq C\bar\rho H\|\sqrt\rho\,\dot u\|_2^2+C\|P\|_4^4.
		\end{aligned}
	\end{equation}
	Taking $p=3$ in \eqref{eq:pressure-power} and  using \eqref{eq:elliptic-Lp},
	\eqref{eq:flux-Lp} and
	$P^3|F|\le\eta P^4+C_\eta|F|^4$,  with $\eta>0$ sufficiently small,
	we obtain
	\begin{equation}\label{eq:pressure-3}
		\frac{\mathrm d}{\mathrm dt}\|P\|_3^3+\frac{3\gamma-1}{2(2\mu+\lambda)}\|P\|_4^4
		\le C\|F\|_4^4\le C\|F\|_2^2\|\nabla F\|_2^2
		\le C\bar\rho H\|\sqrt\rho\,\dot u\|_2^2.
	\end{equation}
	 Adding $K_3$ times \eqref{eq:pressure-3} to  \eqref{eq:accestimate}, we obtain
	\begin{equation}\label{eq:Y-differential}
		\frac{\mathrm d}{\mathrm dt}[\|\sqrt\rho\,\dot u\|_2^2+K_3\|P\|_3^3]
		+c\left(\|\nabla\dot u\|_2^2+\|P\|_4^4\right)
		\le C\bar\rho H\|\sqrt\rho\,\dot u\|_2^2
		\le C\bar\rho H[\|\sqrt\rho\,\dot u\|_2^2+K_3\|P\|_3^3],
	\end{equation}
	provided that
  $$K_3\frac{2\gamma-1}{2(2\mu+\lambda)}\ge C+1.$$
	Multiplying \eqref{eq:Y-differential} by $t$ and using  \eqref{eq:H-decay}, we get
	\begin{equation}\label{eq:deffi-Y}
		\begin{aligned}
			&\frac{\mathrm d}{\mathrm dt}(t[\|\sqrt\rho\,\dot u\|_2^2+K_3\|P\|_3^3])
			+ct\left(\|\nabla\dot u\|_2^2+\|P\|_4^4\right)
			\\&\le\|\sqrt\rho\,\dot u\|_2^2+K_3\|P\|_3^3+C\bar\rho\,tH\|\sqrt\rho\,\dot u\|_2^2\\
			&\le\|\sqrt\rho\,\dot u\|_2^2+K_3\|P\|_3^3+C\bar\rho E_0\|\sqrt\rho\,\dot u\|_2^2
			\\&\le C\|\sqrt\rho\,\dot u\|_2^2+K_3\|P\|_3^3,
		\end{aligned}
	\end{equation}
	where we have used the  boundedness of $\rho E_0$.
	
	Next, we show that $$t[\|\sqrt\rho\,\dot u(t)\|_2^2+K_3\|P(t)\|_3^3]\to0,$$ as $t\to0$.  For almost every
	$s\in[t/2,t]$, \eqref{eq:Y-differential} gives
	\begin{equation}\label{eq:ineqforyt}
		\begin{aligned}\|\sqrt\rho\,\dot u(t)\|_2^2+K_3\|P(t)\|_3^3&\le\|\sqrt\rho\,\dot u(s)\|_2^2+K_3\|P(s)\|_3^3
		\exp\left(C\bar\rho\int_s^tH(\xi)\,\mathrm d\xi\right)
		\\&\le e^{C\bar\rho H_0t}[\|\sqrt\rho\,\dot u(s)\|_2^2+K_3\|P(s)\|_3^3],
\end{aligned}
	\end{equation}
	where we have used \eqref{eq:first-integrated}. Integrating \eqref{eq:ineqforyt} with respect to $s$  over  $ [t/2,t]$, we obtain
	\begin{equation}\label{limit-Y}
	\frac t2[\|\sqrt\rho\,\dot u(t)\|_2^2+K_3\|P(t)\|_3^3]
	\leq e^{C\bar\rho H_0t}\int_{t/2}^t[\|\sqrt\rho\,\dot u(s)\|_2^2+K_3\|P(s)\|_3^3]\,\mathrm ds
	\longrightarrow0, \quad\text{as }t\to0.
	\end{equation}
	We now integrate   \eqref{eq:deffi-Y} over $[\tau, t]$ and let $\tau\to 0^+$. Combining  \eqref{limit-Y} and \eqref{eq:first-integrated}, we obtain  \eqref{eq:Y-t}.
	
	It remains to prove the pointwise decay estimate \eqref{eq:Y-decay}. By  \eqref{eq:H-decay}, we have
	\[
	C\bar\rho\int_{t/2}^tH(\xi)\,\mathrm d\xi
	\leq C\bar\rho E_0
	\log\frac{t+H_0/E_0}{t/2+H_0/E_0}
	\le C\bar\rho E_0\log2\le C.
	\]
	Thus, for almost every $s\in[t/2,t]$,
	\begin{equation}\label{eq:ineqfor-pointwise}
		\begin{aligned}
			\|\sqrt\rho\,\dot u(t)\|_2^2+K_3\|P(t)\|_3^3&\leq [\|\sqrt\rho\,\dot u(s)\|_2^2+K_3\|P(s)\|_3^3]	\exp\left(C\bar\rho\int_{s}^tH(\xi)\,\mathrm d\xi\right)\\&\leq [\|\sqrt\rho\,\dot u(s)\|_2^2+K_3\|P(s)\|_3^3]	\exp\left(C\bar\rho\int_{t/2}^tH(\xi)\,\mathrm d\xi\right)\\&\leq C \|\sqrt\rho\,\dot u(s)\|_2^2+K_3\|P(s)\|_3^3
		\end{aligned}
	\end{equation}
	Integrating \eqref{eq:ineqfor-pointwise} with respect to $s$ over $[t/2,t]$, and using \eqref{eq:first-energy} and \eqref{eq:H-decay}, we have
	\[
	\|\sqrt\rho\,\dot u(t)\|_2^2+K_3\|P(t)\|_3^3\le\frac Ct\int_{t/2}^t[\|\sqrt\rho\,\dot u(s)\|_2^2+K_3\|P(s)\|_3^3]\,\mathrm ds
	\leq\frac CtH(t/2)
	\leq\frac{CE_0}{t(t+H_0/E_0)},
	\]
	which implies \eqref{eq:Y-decay}. This completes the proof.
\end{proof}

\begin{lemma}
	\label{lem:mass-pressure}
	Under the assumption of Proposition \ref{prop:energy}, we have
	\begin{equation}\label{eq:time-integrals}
		\begin{aligned}
			&\int_0^T H(t)\dt
			+\int_0^Tt\left(\|\sqrt\rho\,\dot u\|_2^2+\|P\|_3^3\right)\dt\\
			&\qquad
			+\int_0^Tt^2\left(\|\nabla\dot u\|_2^2+\|P\|_4^4\right)\dt
			\le C(E_0+m_0H_0),
		\end{aligned}
	\end{equation} provided that
	\begin{equation*}\bar\rho E_0 \le \min\left\{1,  \frac{c}{2C}\right\}.\end{equation*}
\end{lemma}
\begin{proof}Similar to the argument of  (3.22) in  \cite{LiXinPre}.
	Noticing that  \eqref{eq:cns}$_2$ gives that
	$$P=(-\Delta)^{-1}\Div(\rho\dot u)+\nu\Div u.$$  Then  using H\"older's inequality and interpolation inequality, we have
 \begin{equation}
		\begin{aligned}
			\|P\|_2^2
			&\leq C\|\rho\dot u\|_{4\gamma/(2\gamma+1)}
			\|P\|_{4\gamma/(4\gamma-1)}
			+C\|\nabla u\|_2\|P\|_2\\
			&\leq C \|\rho\|_{2\gamma}^{1/2}\|\sqrt\rho\,\dot u\|_2\|\rho\|_{4\gamma^2/(4\gamma-1)}^\gamma+C\|\nabla u\|_2\|P\|_2\\
			&\leq C\|\rho\|_{2\gamma}^{1/2}\|\sqrt\rho\,\dot u\|_2 m_0^{1/2}\|\rho\|_{2\gamma}^{\gamma-1/2}+C\|\nabla u\|_2\|P\|_2\\
			&\leq  C\sqrt{m_0}\|\sqrt\rho\,\dot u\|_2\|P\|_2
			+C\|\nabla u\|_2\|P\|_2.
		\end{aligned}
	\end{equation}
	Thus we get
	\begin{equation}\label{eq:mass-pressure}
		\|P\|_2^2\leq  C m_0\|\sqrt\rho\,\dot u\|_2^2
		+C\|\nabla u\|_2^2.
	\end{equation}
	Integrating \eqref{eq:mass-pressure} over $[0,T]$ and using
	\eqref{eq:basic} and \eqref{eq:first-integrated}, we obtain
	\begin{equation}\label{eq:H}
		\int_0^T H(t)\dt\le C(E_0+m_0H_0).
	\end{equation}
	Next, multiplying \eqref{eq:X-coercive} by $tM_0(t)$ where
	$M_0$ is given by (\ref{M}), we have
	\begin{equation}\label{eq:P3estimate}
		\frac{\mathrm d}{\mathrm dt}(tM_0\X)
		+ctM_0\left(\|\sqrt\rho\,\dot u\|_2^2+\|P\|_3^3\right)
		\le M_0\X.
	\end{equation}
	Since $\delta M_0(\delta)\X(\delta)\le C\delta H_0\to0$,
	integrating \eqref{eq:P3estimate} over $[\delta,T]$ and  letting $\delta\to 0^+$  give
	\begin{equation}\label{eq:ineqfor-int-y}
		\int_0^Tt[\|\sqrt\rho\,\dot u(t)\|_2^2+K_3\|P(t)\|_3^3]\dt
		\le C\int_0^T H(t)\dt
		\le C(E_0+m_0H_0).
	\end{equation}
	Finally, multiplying \eqref{eq:Y-differential} by $t^2$ and using \eqref{eq:H-decay}, we obtain
	\begin{equation}\label{eq:ineqfor-ty}
		\begin{aligned}
			&\frac{\mathrm d}{\mathrm dt}(t^2[\|\sqrt\rho\,\dot u(t)\|_2^2+K_3\|P(t)\|_3^3])
			+ct^2\left(\|\nabla\dot u\|_2^2+\|P\|_4^4\right)
			\\&\le2t[\|\sqrt\rho\,\dot u(t)\|_2^2+K_3\|P(t)\|_3^3]+C\bar\rho\,t^2H\|\sqrt\rho\,\dot u\|_2^2\\
			&\le C(1+\bar\rho E_0)t[\|\sqrt\rho\,\dot u(t)\|_2^2+K_3\|P(t)\|_3^3]
			\le Ct[\|\sqrt\rho\,\dot u(t)\|_2^2+K_3\|P(t)\|_3^3],
		\end{aligned}
	\end{equation}
	where we have used the  boundedness of $\rho E_0$. Integrating \eqref{eq:ineqfor-ty} over $(s,t)$ for $s<t$, and letting $s\to 0^+$, and using \eqref{eq:ineqfor-int-y}, we get
	\begin{equation}\label{eq:ineqforpressure}
		\int_0^Tt^2\left(\|\nabla\dot u\|_2^2+\|P\|_4^4\right)\dt
		\le C(E_0+m_0H_0),
	\end{equation}
	where we have used (\ref{limit-Y}). Combining \eqref{eq:H}, \eqref{eq:ineqfor-int-y} and  \eqref{eq:ineqforpressure} gives  \eqref{eq:time-integrals}.
	This completes the proof.
\end{proof}

\begin{proof}[Proof of Proposition~\ref{prop:energy}]
	Using Lemmas~\ref{lem:first-energy}-\ref{lem:mass-pressure},
	we finish the proof of Proposition ~\ref{prop:energy}.
\end{proof}

\subsection{Upper bound of density}\label{sec:density}

\begin{proposition}\label{prop:density}
	Letting $(\rho,u)$ be the strong solution on $[0,T]$, there exists a positive constant $\eps_2$ depending only on $\mu,\lambda,\gamma$ and some other known parameters but independent of T and the initial data, such that if
	\[
	0\le\rho(x,t)\le2\bar\rho,\qquad x\in\R^2,\quad 0\le t\le T,
	\]
	then there holds that
	\[
	0\le\rho(x,t)\le\frac{7}{4}\bar\rho,
	\qquad (x,t)\in\R^2\times[0,T],
	\]
	provided that (\ref{mainassum1}) holds and that
	\begin{equation}\label{mainassum2}
\bar\rho E_0\left(1+(\bar\rho E_0)^3\right)(\bar\rho^\gamma m_0)^{-\frac{6}{5}}  R_*^{\frac{4}{5}}\left[(\bar\rho E_0)^2+  [R_*^4+(\bar\rho^\gamma m_0)^2]  (\bar\rho^{1-\gamma}H_0)^2\right]\le \epsilon_2.
\end{equation}
\end{proposition}

Before we give the proof of Proposition \ref{prop:density}, wes need to prove the following lemmas.


\begin{lemma}\label{lem:moment}
Under the conditions of Proposition  \ref{prop:density}, we have
	\begin{equation}\label{eq:increaseI}
	\begin{aligned}
		\int_{\R^2}\rho(x,t)
		\sqrt{1+\bar\rho^{\gamma+1}|x|^2}\dd
		\le
		m_0+\bar\rho^{\frac{\gamma+1}{2}}
			\int_{\R^2}\rho_0|x|\dd
			+t\sqrt{2\bar\rho^{\gamma+1}m_0E_0}
	\end{aligned}
\end{equation}
and
	\begin{equation}\label{eq:masslowbd}
		\int_{B_{R^*(t)}}\rho(x,t)\dd
		\ge \frac{m_0}{2},
	\end{equation}
	where $R^*(t)$ is given in \eqref{eq:increaseR}.

\end{lemma}

\begin{proof}
	Multiplying the continuity equation by
	$\sqrt{1+\bar\rho^{\gamma+1}|x|^2}$
	and integrating by parts, we obtain
	\begin{equation}\label{eq:It}
		\begin{aligned}
			\frac{\mathrm d}{\mathrm dt}
			\int\rho
			\sqrt{1+\bar\rho^{\gamma+1}|x|^2}\dd
		     &=
			\int\rho u\cdot
			\nabla\sqrt{1+\bar\rho^{\gamma+1}|x|^2}\dd\\
			& \leq
			\bar\rho^{\frac{\gamma+1}{2}}
			\sqrt{m}\|\sqrt\rho\,u\|_2\\
			& \leq
			\sqrt{2\bar\rho^{\gamma+1}m_0E_0},
		\end{aligned}
	\end{equation}
	where we have used (\ref{eq:rho}) and (\ref{eq:basic}).
	
	Integrating \eqref{eq:It} over $(0,t)$ and using
	$\sqrt{1+s^2}\le 1+s$ for $s\ge0$, we obtain
	\begin{equation}\label{eq:moment-growth}
		\begin{aligned}
			\int\rho(x,t)
			\sqrt{1+\bar\rho^{\gamma+1}|x|^2}\dd
			&\le
			\int\rho_0
			\sqrt{1+\bar\rho^{\gamma+1}|x|^2}\dd
			+t\sqrt{2\bar\rho^{\gamma+1}m_0E_0}\\
			&\quad\le
			m_0+\bar\rho^{\frac{\gamma+1}{2}}
			\int\rho_0|x|\dd
			+t\sqrt{2\bar\rho^{\gamma+1}m_0E_0},
		\end{aligned}
	\end{equation} which gives (\ref{eq:increaseI}).
	
	Next, let $R^*(t)>0$ be chosen below. For $|x|\ge R^*(t)$, we have
	\begin{equation}
		\sqrt{1+\bar\rho^{\gamma+1}|x|^2}
		\ge \bar\rho^{(\gamma+1)/2}R^*(t).
	\end{equation}
	Hence, by \eqref{eq:moment-growth},
	\begin{equation}
		\begin{aligned}
			&\int_{|x|\ge R^*(t)}\rho(x,t)\dd\\
			&\quad\le
			\frac{1}{\bar\rho^{\frac{\gamma+1}{2}}R^*(t)}
			\int_{\R^2}\rho(x,t)
			\sqrt{1+\bar\rho^{\gamma+1}|x|^2}\dd\\
			&\quad\le
			\frac{m_0}{\bar\rho^{\frac{\gamma+1}{2}}R^*(t)}
			\left(
			1+\frac{\bar\rho^{\frac{\gamma+1}{2}}}{m_0}
			\int_{\R^2}\rho_0|x|\dd
			+\sqrt{\frac{2\bar\rho^{1-\gamma}E_0}{m_0}}
			\,\bar\rho^\gamma t
			\right)\\
			&\quad=\frac{m_0}{2},
		\end{aligned}
	\end{equation}
	where  we choose
	\begin{equation}\label{eq:increaseR}
		\begin{aligned}
			R^*(t):={}&
			2\bar\rho^{-\frac{\gamma+1}{2}}
			\left(
			1+\frac{\bar\rho^{\frac{\gamma+1}{2}}}{m_0}
			\int_{\R^2}\rho_0|x|\dd
			+\sqrt{\frac{2\bar\rho^{1-\gamma}E_0}{m_0}}
			\bar{\rho}^{\gamma}t\right).
		\end{aligned}
	\end{equation}
	
	By \eqref{eq:rho}, we have
	\begin{equation*}
		\int_{B_{R^*(t)}}\rho(x,t)\dd
		=m_0-\int_{\R^2\setminus B_{R^*(t)}}\rho(x,t)\dd
		\ge\frac{m_0}{2}.
	\end{equation*}
	This proves \eqref{eq:masslowbd} and completes the proof.
\end{proof}

\begin{lemma}\label{lem:bridge}
	Under the assumptions of Proposition~\ref{prop:density},
	the following estimate holds:
	\begin{equation}\label{eq:bridge}
		\begin{aligned}
		\int_{\R^2}\rho|v|^p\dd
		\leq C\bar\rho^{-\frac{\gamma}{5} }
		\left(m_0R\right)^{\frac{4}{5}}\left( m_0^{-1/2}\|\sqrt\rho v\|_2
		+[(\bar\rho^\gamma m_0)^{-1/2}R+1]\|\nabla v\|_2\right)^p,
	\end{aligned}
	\end{equation} for any given $v\in D^{1,2}$ and $\sqrt{\rho}v\in L^2$,
	where
	$$
	R(t)=2\left(
	1+\frac{\bar\rho^{(\gamma+1)/2}}{m_0}
	\int_{\R^2}\rho_0|x|\dd
	+\sqrt{\frac{2\bar\rho^{1-\gamma}E_0}{m_0}}
	\bar{\rho}^{\gamma}t\right).
	$$

\end{lemma}

\begin{proof}
	Fixing $t\in[0,T]$, we denote by $B_{\bar\rho^{-\frac{\gamma+1}{2}}R}$
	the ball centered at the origin with radius
	$\bar\rho^{-\frac{\gamma+1}{2}}R$ where $R=R(t)$.
	Letting $\tilde v$ be the average of $v$ over
	$B_{\bar\rho^{-(\gamma+1)/2}R}$, namely,
	\begin{align}
		\tilde v =\frac{1}{|B_{\bar\rho^{-\frac{\gamma+1}{2}}R}|}\int_{B_{\bar\rho^{-\frac{\gamma+1}{2}}R}}v\dd,
	\end{align}
	and using the a priori assumption $\rho\le2\bar\rho$ in Proposition~\ref{prop:density}, we have
	\begin{equation}\label{eq:rhou}
	\begin{aligned}
		|\tilde v|\left(\int_{B_{\bar\rho^{-\frac{\gamma+1}{2}}R}}\rho\dd\right)^{1/2}&=\|\sqrt{\rho}\tilde v\|_{L^2(B_{\bar\rho^{-\frac{\gamma+1}{2}}R})}\\
		&\leq
		\|\sqrt\rho v\|_2
		+\|\sqrt\rho(v-\tilde v)\|_{L^2(B_{\bar\rho^{-\frac{\gamma+1}{2}}R})} \\
		&\leq C
		\left(\|\sqrt\rho v\|_2
		+\bar\rho^{-\frac{\gamma}{2} }R\|\nabla v\|_2\right),
	\end{aligned}
	\end{equation}
	where the last inequality follows from the following  Poincar\'e's inequality
	$$
	\|v-\tilde v\|_{L^2(B_{ \rho^{-\frac{\gamma+1}{2}}R})}\leq C\rho^{-\frac{\gamma+1}{2}}R \|\nabla v\|_2.
	$$
	Combining \eqref{eq:rhou} with \eqref{eq:masslowbd} and using the fact that $R^*=\bar\rho^{-\frac{\gamma+1}{2}}R$, we get
	\begin{equation}\label{utilde}
		|\tilde v|\leq C m_0^{-1/2}
		\left(\|\sqrt\rho v\|_2
		+\bar\rho^{-\gamma/2}R\|\nabla v\|_2\right).
	\end{equation}
	Using Poincar\'e's inequality and (\ref{utilde}), we obtain
	\[
	\begin{aligned}
		\bar\rho^{\frac{\gamma+1}{2}}
		\|v\|_{L^2(B_{\bar\rho^{-\frac{\gamma+1}{2}}})}&\le
		\bar\rho^{\frac{\gamma+1}{2}}
		\|v-\tilde v\|_{L^2(B_{\bar\rho^{-\frac{\gamma+1}{2}}})}+C|\tilde v|\\
		&\le C\|\nabla v\|_2+C|\tilde v|\\ &\le
C\|\nabla v\|_2+C m_0^{-1/2}
		\left(\|\sqrt\rho v\|_2
		+\bar\rho^{-\gamma/2}R\|\nabla v\|_2\right).
	\end{aligned}
	\]
	Consequently, we have
	\begin{equation}\label{eq:local-anchor}
		\begin{aligned}
			&\bar\rho^{\frac{\gamma+1}{2}}
			\|v\|_{L^2(B_{\bar\rho^{-\frac{\gamma+1}{2}}})}+\|\nabla v\|_2\\&\quad\leq
C m_0^{-1/2}\|\sqrt\rho v\|_2
		+C[(\bar\rho^\gamma m_0)^{-1/2}R+1]\|\nabla v\|_2.
		\end{aligned}
	\end{equation}

	Applying Lemma~\ref{lem:weighted-sobolev} to
	$w(x)=v(\bar\rho^{-(\gamma+1)/2}x)$, changing variables, and using (\ref{eq:local-anchor}), we obtain

		\begin{equation}\label{eq:esitmatev-wei}
			\begin{aligned}
				&\left(
				\int
				\frac{|v(x)|^{5p}}
				{\bigl(e+\bar\rho^{\gamma+1}|x|^2\bigr)
					\bigl[\log\bigl(e+\bar\rho^{\gamma+1}|x|^2\bigr)\bigr]^\theta}
				\,dx
				\right)^{\frac{1}{5p} }
				\\
				&\quad\le
				C_{p,\theta}\bar\rho^{-\frac{\gamma+1}{5p}}
				\left(
				\bar\rho^{\frac{\gamma+1}{2}}
				\|v\|_{L^2(B_{\bar\rho^{-\frac{\gamma+1}{2}}})}
				+\|\nabla v\|_2
				\right)\\&\quad\le C_{p,\theta}\bar\rho^{-\frac{\gamma+1}{5p}}
				\left( m_0^{-1/2}\|\sqrt\rho v\|_2
		+[(\bar\rho^\gamma m_0)^{-1/2}R+1]\|\nabla v\|_2\right),
					\end{aligned}
		\end{equation}
		 for a constant $\theta$ given by Lemma \ref{lem:weighted-sobolev}.

	Next, by using H\"older's inequality and the a priori assumption $\rho\le2\bar\rho$, we have
	\begin{equation}\label{eq:rhou1}
	\begin{aligned}
		\left(\int\rho|v|^p\dd\right)^{\frac{1}{p} }
		&\leq (2\bar\rho)^{\frac{1}{5p} }
		\left( \int\rho(x,t)
		\sqrt{1+\bar\rho^{\gamma+1 }|x|^2} \dd\right)^{\frac{4}{5p}}
		\left(\int |v|^{5p}
		(1+\bar\rho^{\gamma+1}|x|^2)^{-2}\dd\right)^{\frac{1}{5p}}.
	\end{aligned}
	\end{equation}
	To estimate the last term on the right-hand side of (\ref{eq:rhou1}), we use
		$$
		(1+|x|^2)^{-2}
		\le
		\frac{C_\theta}
		{(e+|x|^2)[\log(e+|x|^2)]^\theta}.
		$$
	Combining  \eqref{eq:rhou1} with \eqref{eq:increaseI}
	and \eqref{eq:esitmatev-wei}, we obtain
	\begin{equation*}
		\begin{aligned}
		&\left(\int\rho|v|^p\dd\right)^{\frac{1}{p} }
		\\ \quad &\leq C\bar\rho^{-\frac{\gamma}{5p} }
		\left(m_0+\bar\rho^{\frac{\gamma+1}{2}}
			\int_{\R^2}\rho_0|x|\dd
			+t\sqrt{2\bar\rho^{\gamma+1}m_0E_0}\right)^{\frac{4}{5p}}
				\\ \quad &\quad\times\left( m_0^{-1/2}\|\sqrt\rho v\|_2
		+[(\bar\rho^\gamma m_0)^{-1/2}R+1]\|\nabla v\|_2\right)\\
&\leq C\bar\rho^{-\frac{\gamma}{5p} }
		\left(m_0R\right)^{\frac{4}{5p}}\left( m_0^{-1/2}\|\sqrt\rho v\|_2
		+[(\bar\rho^\gamma m_0)^{-1/2}R+1]\|\nabla v\|_2\right).
	\end{aligned}
	\end{equation*}
	This completes the proof.
\end{proof}

We can now prove Proposition~\ref{prop:density}.
\begin{proof}[Proof of Proposition~\ref{prop:density}]
	As in \cite{Wen},	for any given
	$(x,t)\in\mathbb{R}^2\times[0,T]$, denote
	\[
	\rho^\delta(y,s)
	=
	\rho(y,s)
	+\delta\exp\left\{
	-\int_0^s
	\operatorname{div}u\bigl(X(\tau;x,t),\tau\bigr)
	\,\mathrm d\tau
	\right\}
	>0,
	\]
	where $X(s;x,t)$ is given by
	\[
	\begin{cases}
		\displaystyle
		\frac{\mathrm d}{\mathrm ds}X(s;x,t)
		=
		u\bigl(X(s;x,t),s\bigr),
		& 0\le s<t,\\[6pt]
		X(t;x,t)=x.
	\end{cases}
	\]
	It is easy to verify that
	\[
	\frac{\mathrm d}{\mathrm ds}
	\rho^\delta\bigl(X(s;x,t),s\bigr)
	+
	\rho^\delta\bigl(X(s;x,t),s\bigr)
	\operatorname{div}u\bigl(X(s;x,t),s\bigr)
	=0,
	\]
	due to \eqref{eq:cns}$_1$. This gives
	\begin{equation}\label{eq:trajectory-ode}
		Y'(s)=g(s)+b'(s),
	\end{equation}
	where
	\[
	Y(s)=\ln\rho^\delta\bigl(X(s;x,t),s\bigr),
	\qquad
	g(s)=
	-\frac{\rho^\gamma\bigl(X(s;x,t),s\bigr)}
	{2\mu+\lambda},
	\]
	\[
	b(s)=
	-\frac{1}{2\mu+\lambda}
	\int_0^s F\bigl(X(\tau;x,t),\tau\bigr)
	\,\mathrm d\tau,
	\]
	and $F=(2\mu+\lambda)\operatorname{div}u-\rho^\gamma.$

	For $0\le t_1<t_2\le {\min\{\bar\rho^{-\gamma},T\}}$, we have
	\begin{equation}\label{eq:orbit-to-norm1}
		\begin{aligned}
			b(t_2)-b(t_1)
			&\leq\frac{1}{2\mu+\lambda}\int_{t_1}^{t_2}\|F(s)\|_\infty\,ds\\
			&\leq\frac{1}{2\mu+\lambda}\int_0^{\min\{\bar\rho^{-\gamma},T\}}
			\|F(s)\|_\infty\dt\\&=: I_{1}.
		\end{aligned}
	\end{equation}

	For $I_1$, by Proposition~\ref{prop:energy} and Lemma~\ref{lem:elliptic},
	we have
	\begin{equation}\label{eq:estim1}
	\|F\|_2
	\leq CE_0^{\frac{1}{2}}t^{-\frac{1}{2}},
	\qquad
	\|\nabla F\|_4
	\leq C\|\rho\dot u\|_4.
	\end{equation}
	 Applying \eqref{eq:bridge} with $p=4$ and $v=\dot u$, we obtain
	\begin{equation}\label{eq:bridge-four-growth1}
		\begin{aligned}
			\|\rho\dot u\|_4&\le \bar{\rho}^{\frac{3}{4}}
			\left(\int_{\R^2}\rho|\dot u|^4\dd\right)^{\frac{1}{4} }\\
			&\le C\bar\rho^{\frac{3}{4}-\frac{\gamma}{20} }
			\left(m_0R\right)^{\frac{1}{5}}\left( m_0^{-1/2}\|\sqrt\rho \dot u\|_2
			+[(\bar\rho^\gamma m_0)^{-1/2}R+1]\|\nabla \dot u\|_2\right) \\
			&\le C\bar\rho^{\frac{3-\gamma}{4}} (\bar\rho^\gamma m_0)^{\frac{1}{5}}R^{\frac15}\left((\bar\rho^\gamma m_0)^{-\frac{1}{2}} \bar\rho^{\frac{\gamma}{2}}\|\sqrt\rho\,\dot u\|_2+[(\bar\rho^\gamma m_0)^{-1/2}R+1]\|\nabla \dot u\|_2\right)\\
				&\le C\bar\rho^{\frac{3-\gamma}{4}}  (\bar\rho^\gamma m_0)^{-\frac{3}{10}}R^{\frac15} \bar\rho^{\frac{\gamma}{2}}\|\sqrt\rho\,\dot u\|_2+C\bar\rho^{\frac{3-\gamma}{4}}[(\bar\rho^\gamma m_0)^{-\frac{3}{10}}R+(\bar\rho^\gamma m_0)^{\frac{1}{5}}]R^{\frac15}\|\nabla \dot u\|_2 .
		\end{aligned}
	\end{equation}
	Then using the Gagliardo--Nirenberg
	inequality, \eqref{eq:estim1}, \eqref{eq:bridge-four-growth1},  the  H\"older's inequality and
	Lemmas~\ref{lem:second-energy} and \ref{lem:mass-pressure},
	we obtain
	\begin{equation}
		\begin{aligned}
			I_1
			&\leq C\int_0^{\min\{\bar\rho^{-\gamma},T\}}\|F\|_2^{\frac{1}{3}}
			\|\nabla F\|_4^{\frac{2}{3}}\,\mathrm dt
			\leq C\int_0^{\min\{\bar\rho^{-\gamma},T\}}E_0^{\frac{1}{6}} t^{-\frac{1}{6}} \|\rho\dot u(t)\|_4^{\frac{2}{3}} \,\mathrm dt \\
			&\leq 	C\bar\rho^{\frac{3-\gamma}{6}}
			E_0^{\frac{1}{6}}(\bar\rho^\gamma m_0)^{-\frac{1}{5}}  \bar\rho^{\frac{\gamma}{3}}R_*^{\frac{2}{15}}\int_0^{\min\{\bar\rho^{-\gamma},T\}}t^{-\frac{1}{2}}(t\|\sqrt\rho\,\dot u\|_2^2)^{\frac13} \dd\\
			&+ C\bar\rho^{\frac{3-\gamma}{6}}
			E_0^{\frac{1}{6}} [(\bar\rho^\gamma m_0)^{-\frac{3}{10}}R_*+(\bar\rho^\gamma m_0)^{\frac{1}{5}}]^\frac{2}{3}R_*^{\frac{2}{15}}\int_0^{\min\{\bar\rho^{-\gamma},T\}} t^{-\frac{1}{2}}  (t\|\nabla \dot u\|_2^2)^{\frac13}  \dd\\
			&\le C\bar\rho^{\frac{3-\gamma}{6}}
			E_0^{\frac{1}{6}}(\bar\rho^\gamma m_0)^{-\frac{1}{5}}  \bar\rho^{\frac{\gamma}{3}}R_*^{\frac{2}{15}}  \left(
			\int_0^{\bar\rho^{-\gamma}}t^{-\frac{3}{4}}\,\mathrm dt
			\right)^{\frac{2}{3}} \left(
			\int_0^T
			t
			 \|\sqrt\rho\,\dot u\|_2^2	\dt
			\right)^{\frac{1}{3}}\\
			&+ C\bar\rho^{\frac{3-\gamma}{6}}
			E_0^{\frac{1}{6}} [(\bar\rho^\gamma m_0)^{-\frac{3}{10}}R_*+(\bar\rho^\gamma m_0)^{\frac{1}{5}}]^\frac{2}{3}R_*^{\frac{2}{15}} \left(
			\int_0^{\bar\rho^{-\gamma}}t^{-\frac{3}{4}}\,\mathrm dt
			\right)^{\frac{2}{3}} \left(
			\int_0^T
			t
			\|\nabla\dot u\|_2^2	\dt
			\right)^{\frac{1}{3}}\\
			&\le C(\bar\rho E_0)^{\frac{1}{6}}(\bar\rho^\gamma m_0)^{-\frac{1}{5}}  R_*^{\frac{2}{15}}  \left(\bar\rho E_0+\bar\rho m_0H_0
			\right)^{\frac{1}{3}}\\
			&+ C(\bar\rho E_0)^{\frac{1}{6}} [(\bar\rho^\gamma m_0)^{-\frac{3}{10}}R_*+(\bar\rho^\gamma m_0)^{\frac{1}{5}}]^\frac{2}{3}R_*^{\frac{2}{15}}  (\bar\rho^{1-\gamma}H_0)^{\frac{1}{3}}
			\\
			&\le C(\bar\rho E_0)^{\frac{1}{6}}(\bar\rho^\gamma m_0)^{-\frac{1}{5}}  R_*^{\frac{2}{15}}\left[(\bar\rho E_0)^{\frac{1}{3}}+  [R_*^\frac{2}{3}+(\bar\rho^\gamma m_0)^{\frac{1}{3}}]  (\bar\rho^{1-\gamma}H_0)^{\frac{1}{3}}\right],
		\end{aligned}
	\end{equation} where
\begin{equation}\label{Rstar}
	R_*=2\left(
	1+\frac{\bar\rho^{(\gamma+1)/2}}{m_0}
	\int_{\R^2}\rho_0|x|\dd
	+\sqrt{\frac{2\bar\rho^{1-\gamma}E_0}{m_0}}\right)\ge R(t)
	\end{equation} for all $t\in[0,\min\{\bar\rho^{-\gamma},T\}]$.
Thus taking
\[
N_0=C(\bar\rho E_0)^{\frac{1}{6}}(\bar\rho^\gamma m_0)^{-\frac{1}{5}}  R_*^{\frac{2}{15}}\left[(\bar\rho E_0)^{\frac{1}{3}}+  [R_*^\frac{2}{3}+(\bar\rho^\gamma m_0)^{\frac{1}{3}}]  (\bar\rho^{1-\gamma}H_0)^{\frac{1}{3}}\right],
\quad N_1=0,
\]
and noticing that $g(s)\le 0=N_1$ for all s, Lemma ~\ref{lem:zlotnik} gives
\begin{equation*}
Y(s)\le Y(0)+N_0,
\end{equation*}
which yields
\[
\rho^\delta(x,s)
\le
(\bar \rho+\delta)
\exp\{N_0\}.
\]
Letting $\delta$ tend to zero, we have
\begin{equation}\label{upprho1}
\rho(x,t)\leq\bar \rho\exp\{N_0\}
\end{equation} for all $x\in\mathbb{R}^2$ and $t\in[0,\min\{\bar\rho^{-\gamma},T\}]$.
	Using the smallness assumption \eqref{eq:main-small}, such that
\begin{equation}\label{small1}
N_0 \leq \log\frac{3}{2}.
\end{equation}
 Therefore, we have
\begin{equation}\label{eq:est-rho1}
	\rho(x,t)\leq \frac{3}{2}\bar{\rho}
\end{equation} for all $x\in \mathbb{R}^2$ and $t\in[0,\min\{\bar\rho^{-\gamma},T\}]$.

For ${\min\{\bar\rho^{-\gamma},T\}}\le t_1<t_2\le T$, we have
\begin{equation}\label{eq:orbit-to-norm2}
		\begin{aligned}
			b(t_2)-b(t_1)
			&\leq\frac{1}{2\mu+\lambda}\int_{t_1}^{t_2}\|F(s)\|_\infty\,ds\\
&\leq\frac{\bar\rho^{-\gamma}}{4(2\mu+\lambda)}\int_{t_1}^{t_2}\|F(s)\|^2_\infty\,ds+\frac{\bar\rho^\gamma}{2\mu+\lambda}(t_2-t_1)\\
			&\leq\frac{\bar\rho^{-\gamma}}{4(2\mu+\lambda)}
			\int_{\min\{\bar\rho^{-\gamma},T\}}^T\|F(s)\|_\infty^2\ ds
			+\frac{\bar\rho^\gamma}{2\mu+\lambda}(t_2-t_1)\\
			&=: I_2
			+\frac{\bar\rho^\gamma}{2\mu+\lambda}(t_2-t_1).
		\end{aligned}
	\end{equation}
		For $I_2$,
		 using the a priori assumption $\rho\leq 2\bar{\rho}$ and \eqref{eq:Y-decay}, we have
		\begin{equation}\label{eq:nabla_F}
		\|\nabla F\|_2
		\le C\|\rho\dot u\|_2
		\le C\bar\rho^{\frac{1}{2}}\|\sqrt\rho\,\dot u\|_2
		\le C\bar\rho^{\frac{1}{2}}E_0^{\frac{1}{2}}t^{-1}.
		\end{equation}
		Hence, by the Gagliardo--Nirenberg inequality, (\ref{eq:estim1}) and \eqref{eq:nabla_F}, we get
		\begin{equation}\label{eq:F-eight}
			\|F\|_8
			\le C\|F\|_2^{\frac{1}{4}}
			\|\nabla F\|_2^{\frac{3}{4}}
			\le C\bar\rho^{\frac{3}{8}}
			E_0^{\frac{1}{2}}t^{-\frac{7}{8}}.
		\end{equation}
		Combining \eqref{eq:F-eight} with
		\eqref{eq:bridge-four-growth1} and applying the
		Gagliardo--Nirenberg inequality, and Lemma~\ref{lem:elliptic},  we obtain
		\begin{equation}\label{eq:estforF}
		\begin{aligned}
			\|F\|_\infty^2
			&\le C\|F\|_8^{\frac{4}{3}}
			\|\nabla F\|_4^{\frac{2}{3}}\le C\bar\rho^{\frac{1}{2}}
			E_0^{\frac{2}{3}}t^{-\frac{7}{6}}
			\|\rho \dot u\|_4^{\frac{2}{3}}\\
&\le C\bar\rho^{\frac{1}{2}}
			E_0^{\frac{2}{3}}t^{-\frac{7}{6}}\left[\bar\rho^{\frac{3-\gamma}{2}}  (\bar\rho^\gamma m_0)^{-\frac{3}{5}}R^{\frac{2}{5}} \bar\rho^\gamma\|\sqrt\rho\,\dot u\|_2^2+C\bar\rho^{\frac{3-\gamma}{2}}[(\bar\rho^\gamma m_0)^{-\frac{3}{5}}R^2+(\bar\rho^\gamma m_0)^{\frac{2}{5}}]R^{\frac{2}{5}}\|\nabla \dot u\|_2^2\right]^\frac{1}{3}
\\
&\le C\bar\rho^{\frac{1}{2}}
			E_0^{\frac{2}{3}}t^{-\frac{7}{6}}\left[\bar\rho^{\frac{3-\gamma}{6}}  (\bar\rho^\gamma m_0)^{-\frac{1}{5}}R^{\frac{2}{15}} \bar\rho^\frac{\gamma}{3}\|\sqrt\rho\,\dot u\|_2^\frac{2}{3}+C\bar\rho^{\frac{3-\gamma}{6}}[(\bar\rho^\gamma m_0)^{-\frac{1}{5}}R^\frac{2}{3}+(\bar\rho^\gamma m_0)^{\frac{2}{15}}]R^{\frac{2}{15}}\|\nabla \dot u\|_2^\frac{2}{3}\right]
\\
&\le CE_0^{\frac{2}{3}}\bar\rho^{\frac{6+\gamma}{6}}  (\bar\rho^\gamma m_0)^{-\frac{1}{5}}R(t)^{\frac{2}{15}}t^{-\frac{3}{2}} \big[t\|\sqrt\rho\,\dot u\|_2^2\big]^\frac{1}{3}+CE_0^{\frac{2}{3}}\bar\rho^{\frac{6-\gamma}{6}}(\bar\rho^\gamma m_0)^{-\frac{1}{5}}R(t)^\frac{4}{5}t^{-\frac{3}{2}}\big[t\|\nabla \dot u\|_2^2\big]^{\frac{1}{3}}\\&+CE_0^{\frac{2}{3}}\bar\rho^{\frac{6-\gamma}{6}}(\bar\rho^\gamma m_0)^{\frac{2}{15}}R(t)^{\frac{2}{15}}t^{-\frac{3}{2}}\big[t\|\nabla \dot u\|_2^2\big]^{\frac{1}{3}}.
		\end{aligned}
		\end{equation}
		By the H\"older's inequality and
		Lemmas~\ref{lem:second-energy} and \ref{lem:mass-pressure},
	\eqref{eq:estforF} yields
		\begin{equation}\label{I2}
			\begin{aligned}
			I_2	&\leq C\bar\rho^{-\gamma}
				\int_{\min\{\bar\rho^{-\gamma},T\}}^T\|F(t)\|_\infty^2\dt\leq C\bar\rho^{-\gamma}
				\int_{\bar\rho^{-\gamma}}^\infty\|F(t)\|_\infty^2\dt\\
				&\leq
			CE_0^{\frac{2}{3}}\bar\rho^{\frac{6-5\gamma}{6}}  (\bar\rho^\gamma m_0)^{-\frac{1}{5}}
				\left(\int_{\bar\rho^{-\gamma}}^\infty t^{-\frac{9}{4}}R(t)^{\frac{1}{5}}\dt
				\right)^{\frac{2}{3}}\left(
				\int_0^T
				t\|\sqrt\rho\,\dot u\|_2^2\,\mathrm dt
				\right)^{\frac{1}{3}}\\&
+CE_0^{\frac{2}{3}}\bar\rho^{\frac{6-7\gamma}{6}}(\bar\rho^\gamma m_0)^{-\frac{1}{5}}\left(\int_{\bar\rho^{-\gamma}}^\infty R(t)^\frac{6}{5}t^{-\frac{9}{4}}\dt
				\right)^{\frac{2}{3}}\left(\int_0^Tt\|\nabla \dot u\|_2^2\mathrm dt
				\right)^{\frac{1}{3}}\\&+CE_0^{\frac{2}{3}}\bar\rho^{\frac{6-7\gamma}{6}}(\bar\rho^\gamma m_0)^{\frac{2}{15}}\left(\int_{\bar\rho^{-\gamma}}^\infty t^{-\frac{9}{4}}R(t)^{\frac{1}{5}}\dt\right)^{\frac{2}{3}}\left(\int_0^Tt\|\nabla \dot u\|_2^2\mathrm dt\right)^{\frac{1}{3}}\\
				& \leq C(\bar\rho E_0)^{\frac{2}{3}} (\bar\rho^\gamma m_0)^{-\frac{1}{5}}
				R_*^{\frac{2}{15}}(\bar\rho E_0+\bar\rho m_0H_0)^{\frac{1}{3}}
+C(\bar\rho E_0)^{\frac{2}{3}}(\bar\rho^\gamma m_0)^{-\frac{1}{5}}
				R_*^{\frac{4}{5}}(\bar\rho^{1-\gamma}H_0)^{\frac{1}{3}}\\&+C(\bar\rho E_0)^{\frac{2}{3}}(\bar\rho^\gamma m_0)^{\frac{2}{15}}R_*^{\frac{2}{15}}(\bar\rho^{1-\gamma}H_0)^{\frac{1}{3}}
\\
				& \leq C(\bar\rho E_0)^{\frac{2}{3}}(\bar\rho^\gamma m_0)^{-\frac{1}{5}}R_*^{\frac{2}{15}}\left[(\bar\rho E_0)^\frac{1}{3}
+[(\bar\rho^\gamma m_0)^{\frac{1}{3}}
                 +	R_*^{\frac{2}{3}}](\bar\rho^{1-\gamma}H_0)^{\frac{1}{3}}\right],
			\end{aligned}
		\end{equation}
where we have used (\ref{Rstar}) and that
		\begin{equation*}
			\begin{aligned}
				\int_{\bar\rho^{-\gamma}}^\infty
				t^{-\frac{9}{4}}R^{\frac{6}{5}}(t)\,\mathrm dt&=2^{\frac{6}{5}}\int_{\bar\rho^{-\gamma}}^\infty
				t^{-\frac{9}{4}}\left[
	1+\frac{\bar\rho^{(\gamma+1)/2}}{m_0}
	\int_{\R^2}\rho_0|x|\dd
	+\sqrt{\frac{2\bar\rho^{1-\gamma}E_0}{m_0}}
	\bar{\rho}^{\gamma}t\right]^{\frac{6}{5}}\,\mathrm dt\\
				&=
				2^{\frac{6}{5}}\bar\rho^{\frac{5\gamma}{4}}
				\int_1^\infty s^{-\frac{9}{4}}
				\left(
				1+\frac{\bar\rho^{\frac{\gamma+1}{2}}}{m_0}
				\int_{\R^2}\rho_0|x|\dd
				+\sqrt{\frac{2\bar\rho^{1-\gamma}E_0}{m_0}}\,s
				\right)^{\frac{6}{5}}\,\mathrm ds\\
				&\le
				C\bar\rho^{\frac{5\gamma}{4}}
				\left(
				1+\frac{\bar\rho^{\frac{\gamma+1}{2}}}{m_0}
				\int_{\R^2}\rho_0|x|\dd
				+\sqrt{\frac{2\bar\rho^{1-\gamma}E_0}{m_0}}
				\right)^{\frac{6}{5}}
				\int_1^\infty s^{-\frac{21}{20}}\,\mathrm ds\\
				&\le
				C\bar\rho^{\frac{5\gamma}{4}}
				R_*^{\frac{6}{5}},
			\end{aligned}
		\end{equation*}
and that
	\begin{equation*}
			\begin{aligned}
				\int_{\bar\rho^{-\gamma}}^\infty
				t^{-\frac{9}{4}}R^{\frac{1}{5}}(t)\,\mathrm dt
				\le
				C\bar\rho^{\frac{5\gamma}{4}}
				R_*^{\frac{1}{5}}.
			\end{aligned}
		\end{equation*}

	Thus taking
	\[
	N_0=C(\bar\rho E_0)^{\frac{2}{3}}(\bar\rho^\gamma m_0)^{-\frac{1}{5}}R_*^{\frac{2}{15}}\left[(\bar\rho E_0)^\frac{1}{3}
+[(\bar\rho^\gamma m_0)^{\frac{1}{3}}
                 +	R_*^{\frac{2}{3}}](\bar\rho^{1-\gamma}H_0)^{\frac{1}{3}}\right],
	\qquad N_1=\frac{\bar\rho^\gamma}{2\mu+\lambda},
	\]
and using (\ref{eq:orbit-to-norm2}) and (\ref{I2}), we have
\begin{equation}\label{eq:orbit-to-norm3}
		\begin{aligned}
			b(t_2)-b(t_1)\le N_0+N_1(t_2-t_1)
		\end{aligned}
	\end{equation} for all ${\min\{\bar\rho^{-\gamma},T\}}\le t_1<t_2\le T$.

	Recalling that
	\[
	g=g(Y):=
	-\frac{
		\left(
		e^Y-\delta\exp\left\{
		-\int_0^s
		\operatorname{div}u\bigl(X(\tau;x,t),\tau\bigr)
		\,\mathrm d\tau
		\right\}
		\right)^\gamma
	}{2\mu+\lambda},
	\]
	it is obvious that $g(\infty)=-\infty$. Denoting
	\begin{equation}\label{eq:Y-delta-threshold}
		\bar{Y}_\delta
		:=
		\ln\left\{
		\bar\rho
		+\delta\exp\left\{
		\int_0^T
		\|\operatorname{div}u(\cdot,\tau)\|_{L^\infty}
		\,\mathrm d\tau
		\right\}
		\right\}
	\end{equation}
	and letting $Y\ge\bar{Y}_\delta$, there holds that
	\begin{equation}\label{eq:g-negative-bound}
		\begin{aligned}
			g(Y)
			&=
			-\frac{
				\left(
				e^Y-\delta\exp\left\{
				-\int_0^s
				\operatorname{div}u\bigl(X(\tau;x,t),\tau\bigr)
				\,\mathrm d\tau
				\right\}
				\right)^\gamma
			}{2\mu+\lambda}
			\\[6pt]
			&\le
			-\frac{
				\left(
				e^Y-\delta\exp\left\{
				\int_0^T
				\|\operatorname{div}u(\cdot,\tau)\|_{L^\infty}
				\,\mathrm d\tau
				\right\}
				\right)^\gamma
			}{2\mu+\lambda}
			\\[6pt]
			&\le -N_1.
		\end{aligned}
	\end{equation}
	By virtue of Lemma~\ref{lem:zlotnik} and (\ref{eq:est-rho1}), we have
\begin{equation*}
Y(t)\le \max\{Y({\min\{\bar\rho^{-\gamma},T\}}),\bar{Y}_\delta\}+N_0
\end{equation*}
for all $t\in[{\min\{\bar\rho^{-\gamma},T\}}, T]$,
which yields
	\begin{equation*}
	\rho^\delta(x,t)
	\le
	\max\left\{\frac{3}{2}
	\bar \rho+\delta,\,
	\exp\{\bar{Y}_\delta\}
	\right\}
	\exp\{N_0\}.
	\end{equation*}
	Letting $\delta$ tend to zero, we have
\begin{equation}\label{upprho2}
	\rho(x,t)\leq\frac{3}{2}\bar \rho\exp\{N_0\}
\end{equation}
for all $x\in\mathbb{R}^2$ and $t\in[{\min\{\bar\rho^{-\gamma},T\}}, T]$.

	Using the smallness assumption \eqref{eq:main-small}, such that
\begin{equation}\label{assum2}\frac{3}{2} \exp\{N_0\}\leq  \frac{7}{4},\end{equation}
	we have
   \begin{equation}\label{eq:est-rho2}
	\rho(x,t)\leq \frac{7}{4}\bar{\rho}
	\end{equation} for all $x\in\mathbb{R}^2$ and $t\in[{\min\{\bar\rho^{-\gamma},T\}}, T]$.

\medskip

Based on (\ref{eq:est-rho1}) and (\ref{eq:est-rho2}), we have
\begin{equation}\label{eq:est-rho3}
	\rho(x,t)\leq \frac{7}{4}\bar{\rho}
	\end{equation} for all $x\in\mathbb{R}^2$ and $t\in[0, T]$, provided that
there exists a small positive constant $\epsilon_2$ depending only on $\mu$, $\lambda$, $\gamma$ and some other known parameters but independent of T and the initial data, such that
\begin{equation}\label{N0}
(\bar\rho E_0)^{\frac{1}{6}}\left(1+(\bar\rho E_0)^{\frac{1}{2}}\right)(\bar\rho^\gamma m_0)^{-\frac{1}{5}}  R_*^{\frac{2}{15}}\left[(\bar\rho E_0)^{\frac{1}{3}}+  [R_*^\frac{2}{3}+(\bar\rho^\gamma m_0)^{\frac{1}{3}}]  (\bar\rho^{1-\gamma}H_0)^{\frac{1}{3}}\right]\le \epsilon_2^{\frac{1}{6}}
\end{equation} and that both (\ref{small1}) and (\ref{assum2}) hold.
 The proof of Proposition~\ref{prop:density} is complete.
\end{proof}

\subsection{The higher-order estimates and  the proof of Theorem \ref{thm:main}}\label{sec:continuation}

With Propositions \ref{prop:energy} and \ref{prop:density}, the upper bound of density in (\ref{eq:bootstrap}) holds for all $x\in\R^2$ and $t\in[0, T^*)$. Therefore the corresponding energy estimates in Propositions \ref{prop:energy} and \ref{prop:density} hold on time interval $[0, T^*)$. Consequently, the higher-order estimates and weighted density estimate can  be
proved similarly in  \cite{LiXinPre}, we only present them in the following lemma.

\begin{lemma}\label{lem:finite-regularity}
	Under the assumption of Propositions \ref{prop:energy} and \ref{prop:density}, we have
	\begin{align}
		\sup_{t\le T}\| \rho\|_{L^1\cap H^1\cap W^{1,q}}
		+\int_0^T \bigl(\|\nabla u\|^2_{H^1}+\|\nabla^2u\|_q^{\frac{q+1}{q}}+t\|\nabla^2u\|_q^2\bigr)\dt
		&\le C_T,\label{eq:finite-derivatives} \\
		\sup_{t\le T}\|\bar x^a\rho\|_{L^1\cap H^1\cap W^{1,q}}&\le C_T,  \label{eq:finite-weight}
	\end{align}where $C_T$ may depend on $T<T^*$.
\end{lemma}

\begin{proof}[Proof of Theorem~\ref{thm:main}]
 	Choosing
	\begin{equation}\label{eq:epsilon-choice}
		\varepsilon_*
		=\eps_1+\eps_2,
	\end{equation}where $\eps_1$ and $\eps_2$ are given by
	Propositions~\ref{prop:energy} and~\ref{prop:density}, respectively, (\ref{eq:bootstrap}) and all estimates in
	Proposition \ref{prop:energy} and Lemma \ref{lem:finite-regularity} hold on every interval
	$[0,T]$ for any $T<T^*$. With local-in-time result (Proposition \ref{lem:local}) and these global a priori estimates, it arrives at that $T^*=\infty.$ The proof of Theorem~\ref{thm:main} is complete.
	
\end{proof}

	\section*{Acknowledgements}
The authors would like to thank Professor Jing Li for some kindly and valuable comments. This work was supported by the National Natural Science Foundation of China $\#$ 12471209 and by Guangzhou Basic and Applied
Research Projects \#SL2024A04J01206.   J. Wu is supported by   the National Natural Science Foundation of China  $\#$12671258.

\end{document}